\documentclass[10pt]{amsart}
\usepackage{mathrsfs}
\usepackage{amsfonts} 
\usepackage{graphicx,latexsym,bm,amsmath,amssymb,verbatim,multicol,lscape}
\theoremstyle{definition}

\theoremstyle{remark}

\numberwithin{equation}{section}

\begin{document}
\title[Nair's and Farhi's identities are equivalent]
{Nair's and Farhi's identities involving the least common multiple of binomial coefficients are equivalent}%
\author{Shaofang Hong}
\address{Mathematical College, Sichuan University, Chengdu 610064, P.R. China}
\email{sfhong@scu.edu.cn, s-f.hong@tom.com, hongsf02@yahoo.com }
\thanks{The research was supported partially by Program for New Century
Excellent Talents in University Grant \# NCET-06-0785}

\date{\today}%
\begin{abstract}
In 1982, Nair proved the identity: ${\rm lcm}({k \choose 1}, 2{k
\choose 2}, ..., k{k \choose k}) ={\rm lcm}(1, 2, ..., k),\\ \forall
k\in \mathbb{N}.$ Recently, Farhi proved a new identity: ${\rm
lcm}({k \choose 0}, {k \choose 1}, ..., {k \choose k}) =\frac{{\rm
lcm}(1, 2, ..., k+1)}{k+1},\\ \forall k\in \mathbb{N}.$ In this
note, we show that Nair's and Farhi's identities are equivalent.
\end{abstract}

\maketitle

Throughout this note, let ${\mathbb N}$ denote the set of
nonnegative integers. Define ${\mathbb N}^*:={\mathbb
N}\setminus\{0\}$. There are lots of known results about the least
common multiple of a sequence of positive integers. The most
renowned is nothing else than an equivalent of the prime number
theory; it says that $\log{\rm lcm} (1,2,...,n)\sim n$ as $n$
approaches infinity (see, for instance \cite{[HW]}), where ${\rm
lcm}(1,2,\cdots,n)$ means the least common multiple of $1, 2, ...,
n$. Some authors found effective bounds for lcm$(1,2,...,n)$. Hanson
\cite{[Ha]} got the upper bound lcm$(1,2,...,n)\le 3^n (\forall n\ge
1)$. Nair \cite{[N]} obtained the lower bound ${\rm
lcm}(1,2,\cdots,n)\ge 2^n (\forall n\ge 9)$. Nair \cite{[N]} also
gave a new nice proof for the well-known estimate ${\rm
lcm}(1,2,\cdots,n)\ge 2^{n-1} (\forall n\ge 1)$. Hong and Feng
\cite{[HF]} extended this inequality to the general arithmetic
progression, which confirmed Farhi's conjecture \cite{[F1]}.
Regarding to many other related questions and generalizations of the
above results investigated by several authors, we refer the
interested reader to \cite{[BKS]}, \cite{[FK]},
\cite{[HK]}-\cite{[HY]}.

By exploiting the integral $\int_0^1x^{m-1}(1-x)^{n-m}dx$, Nair
\cite{[N]} showed the following identity involving the binomial
coefficients:\\

{\bf Theorem 1.} (Nair \cite{[N]}) {\it For any $n\in \mathbb{N}^*$,
we have}
$${\rm lcm}({n \choose 1}, 2{n \choose 2}, ..., n{n \choose n})
={\rm lcm}(1, 2, ..., n).$$

Recently, by using Kummer's theorem on the $p$-adic valuation of
binomial coefficients (\cite{[K]}), Farhi \cite{[F2]} provided an
elegant $p$-adic proof to the following new interesting identity
involving the binomial coefficients:\\

{\bf Theorem 2.} (Farhi \cite{[F2]}) {\it For any $n\in \mathbb{N}$,
we have}
$${\rm lcm}({n \choose 0}, {n \choose 1}, ..., {n \choose n})
=\frac{{\rm lcm}(1, 2, ..., n+1)}{n+1}.$$

In this note, we will show that Theorem 1 is equivalent to Theorem
2. Evidently, we can rewrite Theorem 2 as follows:\\

{\bf Theorem 3.} {\it For any $n\in \mathbb{N}^*$, we have}
$$n\cdot {\rm lcm}({{n-1} \choose 0}, {{n-1} \choose 1}, ..., {{n-1} \choose {n-1}})
={\rm lcm}(1, 2, ..., n).$$

Therefore it suffices to show that Theorem 1 is equivalent to Theorem 3.
First, we can easily show the following identity:\\

{\bf Theorem 4.} {\it For any $n\in \mathbb{N}^*$, we have}
$${\rm lcm}({n \choose 1}, 2{n \choose 2}, ..., n{n \choose n})
=n\cdot {\rm lcm}({{n-1} \choose 0}, {{n-1} \choose 1}, ..., {{n-1}
\choose {n-1}}).$$

\begin{proof}
For $1\le t\le n$, since
$$
{n \choose t}=\frac{n!}{t!(n-t)!}=\frac{n\cdot (n-1)!}{t\cdot
(t-1)!(n-t)!}=\frac{n}{t}\cdot
\frac{(n-1)!}{(t-1)!(n-t)!}=\frac{n}{t}\cdot {n-1 \choose t-1},
$$
we infer that
$$
t\cdot {n \choose t}=n\cdot {n-1 \choose t-1}.
$$
It follows immediately that
\begin{align*}
{\rm lcm}({n \choose 1}, 2{n \choose 2}, ..., n{n \choose n})&={\rm
lcm}(n\cdot {n-1 \choose 0}, n\cdot {n-1 \choose 1}, ..., n\cdot
{n-1 \choose n-1})\\
&=n\cdot {\rm lcm}({{n-1} \choose 0}, {{n-1} \choose 1}, ..., {{n-1}
\choose {n-1}})\nonumber
\end{align*}
as required. This completes the proof of Theorem 4.
\end{proof}

From Theorem 4 the equivalence of Theorems 1 and 3 follows
immediately. Finally, noting that ${n\choose {n-t}}={n\choose t}$,
we can further rewrite Theorem 2 as follows.\\

{\bf Theorem 5.} {\it For any $n\in \mathbb{N}^*$, we have
$$n\cdot {\rm lcm}({{n-1} \choose 0}, {{n-1} \choose 1}, ..., {{n-1}
\choose {\lfloor\frac{n-1}{2}}\rfloor}) ={\rm lcm}(1, 2, ..., n),$$
where $\lfloor\frac{n-1}{2}\rfloor$ means the largest integer no
more than $\frac{n-1}{2}$.}

\end{document}